# A Trust Region Method for Pareto Front Approximation

Kwang-Hui Ju[†*], Ju-Song Kim[†]

[†] Faculty of Mathematics, **Kim Il Sung** University,
Pyongyang, Democratic People's Republic of Korea.
[*]Corresponding author Email: kh.ju0117@ryongnamsan.edu.kp

In this paper, we consider black-box multiobjective optimization problems in which all objective functions are not given analytically. In multiobjective optimization, it is important to produce a set of uniformly distributed discrete solutions over the Pareto front to build a good approximation. In black-box biobjective optimization, one can evaluate distances between solutions with the ordering property of the Pareto front. These distances allow to be able to evaluate distribution of all solution points, so it is not difficult to maintain uniformity of solutions distribution. However, problems with more than two objectives do not have ordering property, so it is noted that these problems require other techniques to measure the coverage and maintain uniformity of solutions distribution. In this paper, we propose an algorithm based on a trust region method for the Pareto front approximation in black-box multiobjective optimization. In the algorithm, we select a reference point in the set of non-dominated points by employing the density function and explore area around this point. This ensures uniformity of solutions distribution even for problems with more than two objective functions. We also prove that the iteration points of the algorithm converge to Pareto critical points. Finally, we present numerical results suggesting that the algorithm generates the set of well-distributed solutions approximating the Pareto front, even in the case for the problems with three objective functions.

## 1. Introduction

Multiobjective optimization problems arise in various fields such as engineering, health care, economics and finance as problems that optimize various conflicting objectives (Stewart, 2008; Raunder, M S; Gutjahr, W J; Heidenberger, K; Wagner, J; Pasia, J;, 2010; Tapia, M; Coello, C;, 2007).

In multiobjective optimization, a set of best trade-off solutions should be determined, because there is not optimal solution fitting all objectives in general. These trade-off solutions are defined as Pareto optimal solutions or non-dominated solutions, and the image of the set of those solutions is called the Pareto front. Many iterative solution methods try to find a set of discrete solutions that approximate the Pareto front, and it is necessary that the solutions are uniformly distributed over the Pareto front. In this paper, we focus on ensuring uniformity of solutions distribution in black-box multiobjective optimization.

In the biobjective mesh adaptive direct search algorithm (BIMADS) ( (Audet, C; Savard, G; Zghal, W;, 2008)), the diversity of solutions could be evaluated by computing the Euclidean distance between neighboring solutions based on the ordering property of the Pareto front that biobjective optimization problems possess, and the sparsely distributed area is visited to ensure uniformity of solutions distribution. However, it was noted that the problems with more than three objective functions do not have the ordering property of the Pareto front, so the above approach cannot be used to evaluate uniformity of solutions distribution. This motivated us to study on ensuring uniformity of solutions distribution over the Pareto front for multiobjective optimization problems with more than two objectives.

In black-box optimization, derivative information is not available and a large number of function values are required to obtain viable approximate derivative functions. A common approach for derivative-free methods in multiobjective optimization is direct search such as BIMADS (Audet, C; Savard, G; Zghal, W;, 2008) for biobjective and direct multi-search (DMS) method (Custodio, A; Madeira, J; Vaz, A; Vicente, L;, 2011) for multiobjective. In (Custodio, A; Madeira, J; Vaz, A; Vicente, L;, 2011), they proposed a metric to evaluate the sparsity of the set of solutions to produce uniformly distributed solutions. Their metric is easy to calculate and effective for problems with small number of objectives. However, if the number of objectives increases, the metric accuracy decreases. Besides, if the number of variables increases, the performance of the direct search method deteriorates (Thomann, J; Eichfelder, G;, 2019).

Trust region method can solve the derivative-free optimization problems. There are multiobjective realizations of trust region method. A multiobjective optimization algorithm based on trust region method (TRMP) is proposed in (Villacorta, K; Oliveira, P; Soubeyran, A;, 2014). However, they use the derivative information, so this algorithm is not suitable for black-box problems. A trust region algorithm for heterogeneous multiobjective optimization (TRAHM) in which some of objectives are black-box and others are given analytically, is proposed in (Thomann, J; Eichfelder, G;, 2019). They proved that the algorithm could generate the points which converge to a Pareto critical point. However, TRMP and TRAHM are all try to obtain one solution, not a set of Pareto solutions that approximate the Pareto front. So, they did not mention uniformity of solutions distribution. In this paper, we study on approximating the Pareto front of a multiobjective optimization problems using the trust region method.

Various techniques to evaluate the distribution and ensure uniformity of solutions distribution have been proposed in literatures. The notion of crowding distance for a NSGA-II algorithm is presented in (Deb, 2001), but it is not applicable for problems with more than two objectives. An adaptive grid evolutionary algorithm to maintain diversity is proposed in (Yue, X; Guo, Z; Yin, Y; Liu, X;, 2016). And an approach that converting multiobjective problems into biobjective problems related with diversity and proximity is shown in (Li, M; Yang, S; Liu, X;, 2015). Evolutionary algorithms can produce and improve many approximate solutions in each run, but they require so many function evaluations. So the evolutionary algorithms are not suitable for black-box multiobjective problems that suffer from the limitation of number of function values. Diversity maximization approach (DMA) to ensure the uniformity of solutions distribution is proposed in (Masin, M; Bukchin, Y;, 2008). This approach is capable of solving the mixed-integer and combinational problems. An algorithm for multiobjective optimization based on Hypervolume indicator is proposed in (Bader, J; Zitzler, E;, 2011). In (Wang, H; Jin, Y; Yao, X;, 2017) they focused on addressing the difficulties with existing diversity metrics such as Hypervolume indicator and proposed a new diversity metric in multiobjective optimization inspired by a measure for biodiversity, but it is difficult to measure the density of solutions. (Farhang-Mehr, A; Azarm, S;, 2002) used the notion of influence and density function to define and formulate entropy as a new metric that can quantitatively asses the distribution quality of a solution set in an evolutionary multiobjective optimization algorithm. The computational amount of this metric increases linearly with the number of solutions.

This paper makes the following contributions: (1) We use the density function mentioned in (Farhang-Mehr, A; Azarm, S;, 2002) to select a reference point in the set of non-dominated points and evaluate uniformity of solutions distribution, even in the case for the problems with more than two objective functions. Then, using this density function, we propose a trust region algorithm to approximate the Pareto front of black-box multiobjective optimization problem. We extend the trust region method to a multiobjective optimization as in (Villacorta, K; Oliveira, P; Soubeyran, A;, 2014; Thomann, J; Eichfelder, G;, 2019). But, in contrast to previous algorithms, we generate not one solution but a set of approximate solutions that are uniformly distributed over the entire Pareto front. (2) We prove that the solution set which is produced by the proposed algorithm converges to a set of Pareto critical points under some regular conditions. We use the same assumptions and strategy to prove convergence of our algorithm as (Thomann, J; Eichfelder, G;, 2019). (3) We provide the experimental results to compare to another algorithm in terms of biobjective optimization and show the effectiveness of our algorithm for three objective problems.

The paper is organized as follows.

We introduce the preliminary notions for multiobjective optimization and describe the density function. In the next section, we propose an algorithm to select a reference point among the set of non-dominated solutions. We propose a trust region algorithm for black-box multiobjective optimization using the density function in section 3, followed by convergence analysis of the algorithm in section 4. Numerical results are given in section 5 to show the effectiveness of this algorithm. Some discussion and conclusion follow in section 6.

## 2. Multiobjective optimization

### 2.1. *Problem statement*

In this paper, we consider the/an optimization problem such that

$$\min_{x \in \mathbf{R}^n} F(x) \tag{2.1}$$

$F(x) = (f_1(x), f_2(x), \ldots, f_p(x))^T$ . Here, we assume that the objective functions $f_i : \mathbf{R}^n \rightarrow \mathbf{R}$ are twice continuous differentiable black-box functions and bounded from below for all $i = 1, \ldots, p$ . We use the optimality concept for the above problem (Abraham, A; Jain, L; Goldberg, R;, 2005; Ehrgott, 2005; Audet, C; Savard, G; Zghal, W;, 2008; Thomann, J; Eichfelder, G;, 2019).

**Definition 2.1.** Let $x, y \in \mathbf{R}^n$ be two decision vectors. Then we define $x \preceq y$ ( $x$ weakly dominates $y$ ) if and only if $f_i(x) \leq f_i(y)$ for all $i \in \{1, 2, \ldots, p\}$ , and $x \prec y$ ( $x$ dominates $y$ ) if and only if $x \preceq y$ and $f_j(x) < f_j(y)$ for at least one $j \in \{1, 2, \ldots, p\}$.

**Definition 2.2.** A point $\bar{x} \in \mathbf{R}^n$ is said to be globally (weak) Pareto optimal for (2.1) if and only if there is no $x \in \mathbf{R}^n$ such that $x \prec \bar{x}$ ( $x \preceq \bar{x}$ ). The value $F(\bar{x})$ is called globally (weak) efficient for (2.1) and the image of the set of globally efficient points is called the global (weak) Pareto front or simply the (weak) Pareto front for (2.1).

**Definition 2.3.** A point $\bar{x} \in \mathbf{R}^n$ is said to be locally (weak) Pareto optimal for (2.1) if and only if there exists an open neighborhood of $\bar{x}$, $B(\bar{x})$, such that there is no $x \in B(\bar{x})$ satisfying $x \prec \bar{x}$ ( $x \preceq \bar{x}$ ). And the image of the set of locally (weak) efficient points is called the local (weak) Pareto front.

**Definition 2.4.** Let $F = (f_1, \ldots, f_p)$ be totally differentiable at a point $\bar{x} \in \mathbf{R}^n$. This point is said to be Pareto critical for (2.1) if and only if for every vector $d \in \mathbf{R}^n$ there exists an index $j \in \{1, \ldots, p\}$ such that $\nabla f_j(\bar{x})^T d \geq 0$ holds.

**Lemma 2.5.** ( (Thomann, J; Eichfelder, G;, 2019)). Let $\bar{x} \in \mathbf{R}^n$ be a locally weak efficient for (2.1). Then this point is a Pareto critical for (2.1).

**Lemma 2.6.** ( (Thomann, J; Eichfelder, G;, 2019; Drummond & Svaiter, 2005)). Let $f_i : \mathbf{R}^n \to \mathbf{R}$ be continuously differentiable functions for all $i = 1, \ldots, p$. For the function

$$\omega(x) := -\min_{\|d\| \leq 1} \max_{i=1, \ldots, p} \nabla f_i(x)^T d \qquad (2.2)$$

(i)   The mapping $x \mapsto \omega(x)$ is continuous.

(ii)  It holds $\omega(x) \geq 0$ for all $x \in \mathbf{R}^n$.

(iii) $\omega(x) = 0$ if and only if $x \in \mathbf{R}^n$ is a Pareto critical for (2.1).

**Lemma 2.7.** ( (Thomann, J; Eichfelder, G;, 2019)). Let $x \in \mathbf{R}^n$ be an arbitrary, but fixed point and $d_\omega$ be an optimal solution for (2.2).

(i)   If $x$ is not Pareto critical for (2.1), then $d_\omega$ is a descent direction for (2.1) at the point, i.e. there is a scalar $t_0 > 0$ such that $f_i(x + t d_\omega) < f_i(x)$ for all $t \in (0, t_0]$ and for all $i \in \{1, \ldots, p\}$.

(ii)  There exist scalars $\alpha_i \in [0, 1]$ for $i \in \{1, \ldots, p\}$ with $\sum_{i=1}^{p} \alpha_i = 1$ and $\mu \geq 0$ such that

$$d_\omega = -\mu \sum_{i=1}^{p} \alpha_i \nabla f_i(x)$$

holds. If $x$ is not a Pareto critical for (2.1), then $\|d_\omega\| = 1$. If $x$ is a Pareto critical for (2.1),

$$d_\omega = \sum_{i=1}^{p} \alpha_i \nabla f_i(x) = 0 \, .$$

Furthermore, it holds

$$\omega(x) \le \left\| \sum_{i=1}^{p} \alpha_i \nabla f_i(x) \right\| \, .$$

## 2.2. *Density function*

**Definition 2.8.** Let $\mathbf{i}_1, \mathbf{i}_2, \ldots, \mathbf{i}_p$ be the Cartesian unit vectors along $f_1, f_2, \ldots, f_p$. $G = (\min f_1, \min f_2, \ldots, \min f_p)$ and $B = (\max f_1, \max f_2, \ldots, \max f_p)$ are each called the utopia point and anti-utopia point. $\mathbf{j}_1$ be the unit vector along $GB$. The plane that passes through the origin $G$ and is normal to $\mathbf{j}_1$ is called projection hyper plane for (2.1) and denote $\mathbf{\Pi}_{\mathbf{j}_1}$.

We use the Gram-Schmidt orthogonalization to construct a new unit vectors in $p-1$-dimensional space, projection hyper plane as follows.

$$\mathbf{j}_2 = \frac{\mathbf{i}_2 - (\mathbf{i}_2 \cdot \mathbf{j}_1)\mathbf{j}_1}{\|\mathbf{i}_2 - (\mathbf{i}_2 \cdot \mathbf{j}_1)\mathbf{j}_1\|} \, , \quad \mathbf{j}_3 = \frac{\mathbf{i}_3 - (\mathbf{i}_3 \cdot \mathbf{j}_1)\mathbf{j}_1 - (\mathbf{i}_3 \cdot \mathbf{j}_2)\mathbf{j}_2}{\|\mathbf{i}_3 - (\mathbf{i}_3 \cdot \mathbf{j}_1)\mathbf{j}_1 - (\mathbf{i}_3 \cdot \mathbf{j}_2)\mathbf{j}_2\|} \, , \ldots,$$

$$\mathbf{j}_p = \frac{\mathbf{i}_p - (\mathbf{i}_p \cdot \mathbf{j}_1)\mathbf{j}_1 - (\mathbf{i}_p \cdot \mathbf{j}_2)\mathbf{j}_2 - \cdots - (\mathbf{i}_p \cdot \mathbf{j}_{p-1})\mathbf{j}_{p-1}}{\|\mathbf{i}_p - (\mathbf{i}_p \cdot \mathbf{j}_1)\mathbf{j}_1 - (\mathbf{i}_p \cdot \mathbf{j}_2)\mathbf{j}_2 - \cdots - (\mathbf{i}_p \cdot \mathbf{j}_{p-1})\mathbf{j}_{p-1}\|}$$

The projection of solution $F(x) = (f_1(x), f_2(x), \ldots, f_p(x))$ is defined as

$$\Pr(x) = (y_1, y_2, \ldots, y_{p-1}) = ((F(x) \cdot \mathbf{j}_2), (F(x) \cdot \mathbf{j}_3), \ldots, (F(x) \cdot \mathbf{j}_p)) \, .$$

**Definition 2.9.** ( (Farhang-Mehr, A; Azarm, S;, 2002)). Let $\phi : [0, +\infty) \to \mathbf{R}$ be a decreasing function and $d : \mathbf{\Pi}_{\mathbf{j}_1} \times \mathbf{\Pi}_{\mathbf{j}_1} \to \mathbf{R}_+$ be a distance function. And $X = \{x_1, x_2, \ldots, x_N\}$ is an approximated set of Pareto optimal points for (2.1). Then a function $\Omega_i(y) = \phi(d(\Pr(F(x_i)), \Pr(y)))$, $\Omega_i : \mathbf{R}^p \to \mathbf{R}$, is called influence function of the $i$-th solution point, $x_i$.

**Definition 2.10.** ( (Farhang-Mehr, A; Azarm, S;, 2002)). Let $X = \{x_1, x_2, \ldots, x_N\}$ be an approximated set of Pareto optimal points for (2.1). Then a function $D : \mathbf{R}^p \rightarrow \mathbf{R}$ such that

$$D(y) = \sum_{i=1}^{N} \Omega_i(y) , \quad y \in \mathbf{\Pi}_{\mathbf{j}_1}$$

is called density function for the approximated set of Pareto optimal points.

### 2.2.1. *Decreasing functions.*

Various types of influence functions are constructed followed by the decreasing function $\phi$.

      1. Sharing function

$$\phi(x) = \begin{cases} 1 - (x/\sigma)^{\alpha} , & x \leq \sigma \\ 0, & x > \sigma \end{cases} , \quad \sigma > 0 , \quad \alpha \in \mathbf{N}$$

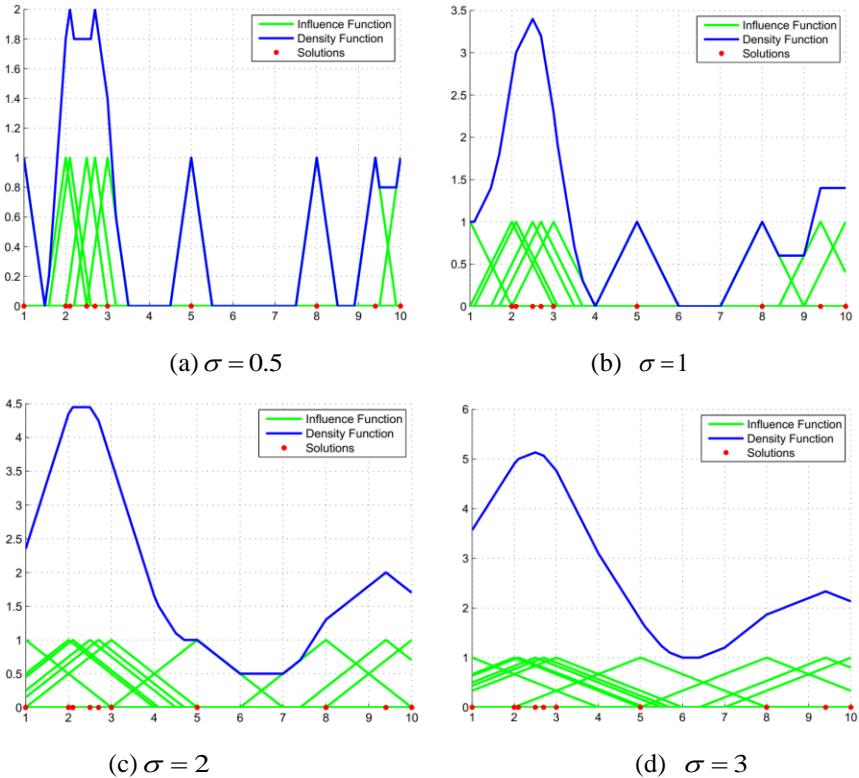

      (a) $\sigma = 0.5$                 (b) $\sigma = 1$

      (c) $\sigma = 2$                 (d) $\sigma = 3$

Fig 1. Influence functions and density functions using sharing function with $\alpha = 1$.

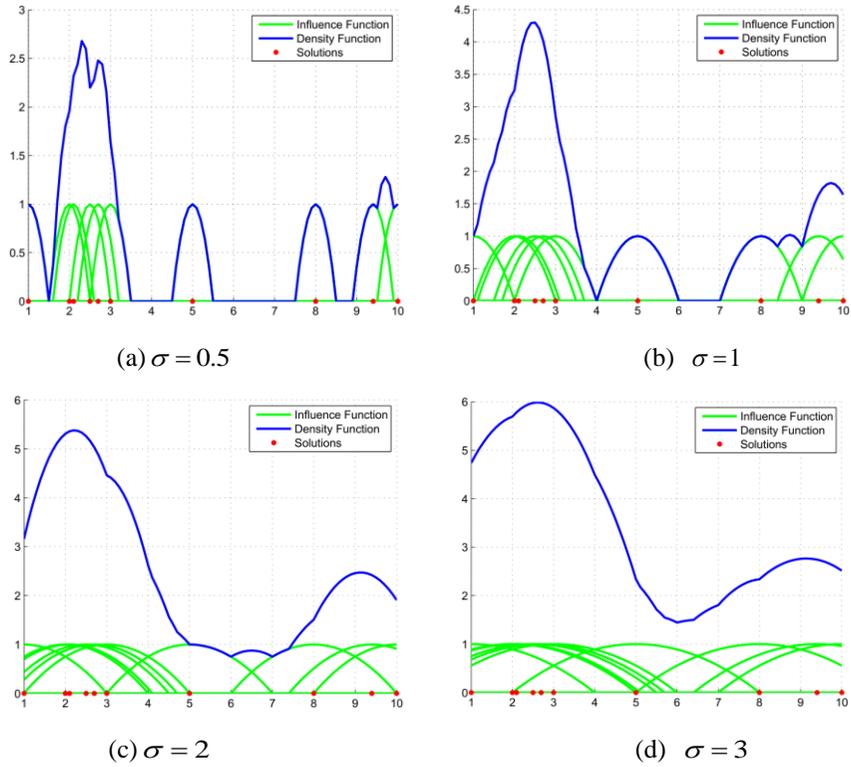

(a) $\sigma = 0.5$          (b) $\sigma = 1$

(c) $\sigma = 2$          (d) $\sigma = 3$

Fig 2. Influence functions and density functions using sharing function with $\alpha = 2$.

2.  Gaussian curve

$$\phi(x) = \exp(-\frac{x^2}{2\sigma^2}), \quad \sigma > 0, \quad x \geq 0$$

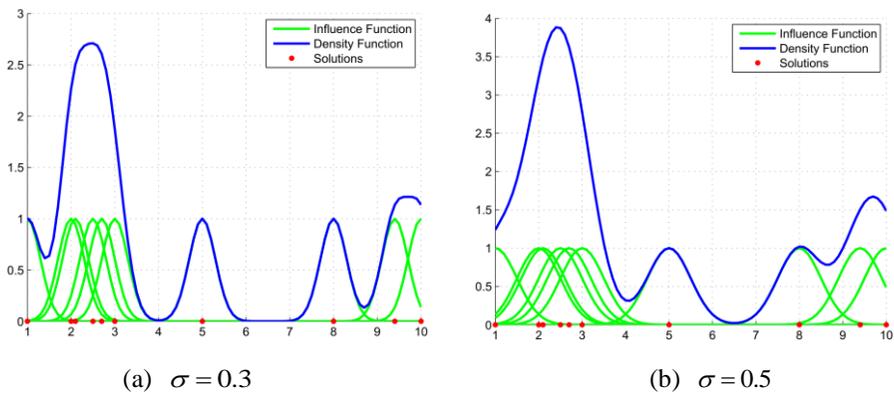

(a)  $\sigma = 0.3$          (b)  $\sigma = 0.5$

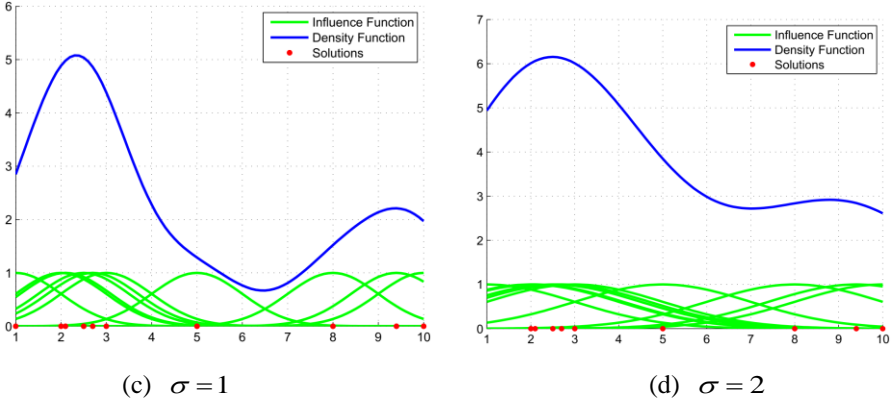

(c) $\sigma = 1$         (d) $\sigma = 2$

Fig 3. Influence functions and density functions using Gaussian curve.

### 2.2.2. *Selecting the reference point*

Let $X = \{x_i, i = 1, 2, \ldots, N\}$ be an approximated set of Pareto optimal points for (2.1). Then the reference point $x_c$ is selected as the following algorithm.

**ALGORITHM 1. (ALGORITHM TO SELECT THE REFERENCE POINT)**

  **Step 1.** Construct influence functions of $x_i \in X$ for all $i = 1, 2, \ldots, N$, $\Omega_i(y)$.

  **Step 2.** Construct a density function for $X = \{x_i, i = 1, 2, \ldots, N\}$,

$D(y) = \sum_{i=1}^{N} \Omega_i(y)$. Calculate the values of density function on $\Pr(F(x_i))$, $\gamma(x_i) = D(\Pr(F(x_i)))$ for all $i = 1, 2, \ldots, N$.

  **Step 3.** $x_c = \arg\min\{\gamma(x_i), i = 1, 2, \ldots, N\}$ is the reference point in $X$.

The reference point obtained by Algorithm 1, is the center point of the most sparsely distributed region in the objective space among $X = \{x_i, i = 1, 2, \ldots, N\}$.

## 3. Trust region method for multiobjective optimization.

In literature, trust region method for single objective optimization are mentioned (Conn, A; Gould, N; Toint, P;, 2000; Conn, A R; Scheinberg, K; Vicente, L N;, 2009). It approximates the objective function by suitable model functions in every iteration. In this method, the computations and models are restricted to a local area in every iteration. We call this area as a trust region. The trust region is defined by

$$B^k := B(x^k, \tilde{\delta}_k) = \{x \in \mathbf{R}^n \mid \|x - x^k\| \leqq \tilde{\delta}_k\}.$$

Here, $x^k$ is the current iteration point, $\tilde{\delta}_k > 0$ is trust region radius and $\|\cdot\| = \|\cdot\|_2$ is the Euclidean norm.

In this paper, we replace each objective function $f_i$ with $i \in \{1, \ldots, p\}$ by a suitable quadratic model $m_i^k : \mathbf{R}^n \to \mathbf{R}$ which satisfies $f_i(x^k) = m_i^k(x^k)$.

So we consider a surrogate problem for (2.1),

$$\min_{x \in \mathbf{R}^n} m^k(x) \tag{3.1}$$

in every iteration $k \in \mathbf{N}$. Here, $m^k(x) = (m_1^k(x), m_2^k(x), \ldots, m_p^k(x))^T$.

## 3.1. *Trial point*

The Pascoletti-Serafini scalarization (Pascoletti, A; Serafini, P;, 1984) is given by

$$\begin{aligned}
&\min && t \\
&\text{s.t.} && F(x^k) - m^k(x) + tr^k \geq 0 \\
& && t \in \mathbf{R} \\
& && x \in B^k
\end{aligned} \tag{3.2}$$

Here, $r^k = (r_1^k, r_2^k, \ldots, r_p^k)^T$, $r_i^k = f_i(x^k) - \min_{x \in B^k} m_i^k(x) \geq 0$ $(i = 1, 2, \ldots, p)$.

**Remark 3.1.** Let's suppose that $x^k$ is not a Pareto-critical for (3.1). Then according to Lemma **2.5**, this point is not a locally weak efficient for (3.1). So, $x^k$ cannot be minimum of any of $m_i^k$, $i \in \{1, \ldots, p\}$. Hence, for every $i \in \{1, \ldots, p\}$, it holds $r_i^k = f_i(x^k) - \min_{x \in B^k} m_i^k(x) > 0$.

**Lemma 3.2.** (Thomann, J; Eichfelder, G;, 2019; EichfelderG., 2008)

a) Let $(t^{k+}, x^{k+})$ be a minimal solution of (3.2). Then $x^{k+}$ is a weakly efficient for $\min_{x \in B^k} m^k(x)$.

b) Let $(t^{k+}, x^{k+})$ be a local minimal solution of (3.2). Then $x^{k+}$ is locally weak efficient for $\min_{x \in B^k} m^k(x)$.

c) Let $x^{k+}$ be a weakly efficient solution for $\min_{x \in B^k} m^k(x)$ and $r^k \in \text{int } \mathbf{R}_+^p$. Then $(0, x^{k+})$ is a minimal solution of (3.2).

**Lemma 3.3.** (Thomann, J; Eichfelder, G;, 2019) If $x^k$ is not a weakly efficient solution for $\min_{x \in B^k} m^k(x)$, then it holds $t^{k+} \in [-1, 0)$ for every minimal solution $(t^{k+}, x^{k+})$ of (3.2).

### 3.2. *Algorithm*

In this paper, we replace every black-box objective functions by quadratic model as in basic trust region method. And we minimize the models with quadratic programming (Thomann, J; Eichfelder, G;, 2019; Conn, A R; Scheinberg, K; Vicente, L N;, 2009).

Since the objective functions are black-box, derivative information cannot be used (Villacorta, K; Oliveira, P; Soubeyran, A;, 2014; Thomann, J; Eichfelder, G;, 2019). So the models like Taylor model cannot be used. To obtain a quadratic model, Interpolation based on the quadratic Lagrangian polynomials is used. For the interpolation of the black-box functions, a set of interpolation points is used. We chose the interpolation points not randomly within the trust region, but by computing with quality criterion called well poisedness described in (Conn, A R; Scheinberg, K; Vicente, L N;, 2009) to be satisfied. In fact, we used the algorithm 6.3 in (Conn, A R; Scheinberg, K; Vicente, L N;, 2009) to get interpolation points.

A quadratic model $m_i(x)$ for $f_i$, $i \in \{1,...,p\}$, is able to be written as

$$f_i(x) \approx m_i(x) = m_i(x_c) + \nabla m_i(x_c)^T (x - x_c) + \frac{1}{2}(x - x_c)^T \nabla^2 m_i(x_c)(x - x_c).$$

We define

$$\omega_m^k(x) := -\min_{\|d\| \le 1} \max_{i=1,...,p} \nabla m_i^k(x)^T d$$

for the model functions.

## ALGORITHM 2. (D-MOTR) (MULTIOBJECTIVE ALGORITHM WITH TRUST REGION METHOD WITH DENSITY FUNCTION)

**Input:** initial point $x^0$, initial trust region radius $\delta_0(x^0)$, lower bound of trust region

$\delta_{tol} > 0$, values for the parameters $0 < \eta_1 \le \eta_2 < 1$, $0 < \gamma_0 < 1$, $\gamma_1$, $\gamma_2$

$(0 < \gamma_1 \le \gamma_2 < 1)$, set $X^0 = \{x^0\}$, $k = 1$.

**Step 1: (Reference point)**

Select the reference point $x^k$ from $X^{k-1}$ using the Algorithm 1. If $x^k = x^{k-1}$

then set $\delta_{k-1}(x^k) = \delta_{k-1}(x^k) \cdot \gamma_0$. Then set $\tilde{\delta}_k = \delta_{k-1}(x^k)$, $j = 1$, $Y^k = \{\varnothing\}$.

**Step 2: (Quadratic models)**

Determine the trust region $B^k = B(x^k, \tilde{\delta}_k) = \{x \in R^n \mid \|x - x^k\| \le \tilde{\delta}_k\}$. Compute the set of sample points $Y^{k,j} = \{y_1^{k,j}, y_2^{k,j}, ..., y_q^{k,j}\} \subset B^k$ according to algorithm 6.3

in (Conn, A R; Scheinberg, K; Vicente, L N;, 2009). Construct quadratic models

$m_i^k(x)$ of $f_i(x)$ for all $i \in \{1, ..., p\}$. We set $Y^k = Y^k \bigcup Y^{k,j}$.

–   If $\tilde{\delta}_k \le \max\{\delta_{tol}, \omega_m^k(x^k)\}$ then go to step 3.

–   Else set $\tilde{\delta}_k = \tilde{\delta}_k \cdot \gamma_0$ and go to step 2.

**Step 3: (Ideal Point)**

   Compute   $s^k = (s_1^k, \ldots, s_p^k)^T$,   $s_i^k = \min_{x \in B^k} m_i^k(x)$,   $i \in \{1, \ldots, p\}$.

**Step 4: (Trial Point)**

   Compute the optimal solution of   (3.2),   $(t^{k+}, x^{k+})^T$.

**Step 5: (Reduction ratio)**

   If there exists an index   $i \in \{1, \ldots, p\}$   such that   $m_i^k(x^k) = m_i^k(x^{k+})$   or   $t^{k+} = 0$,

   then set   $\rho^k = 0$.

   Otherwise compute

$$\rho_i^k = \frac{f_i(x^k) - f_i(x^{k+})}{m_i^k(x^k) - m_i^k(x^{k+})}, \quad i \in \{1, \ldots, p\}.$$

   And set   $\rho^k = \max_{i=1,\ldots,p} \rho_i^k$.

**Step 6: (Updating trust region)**

   For   $x \in Y^k \bigcup \{x^{k+}\}$, set

$$\delta_k(x) \in \begin{cases} [\gamma_1 \widetilde{\delta}_k, \ \gamma_2 \widetilde{\delta}_k], & \rho^k < \eta_1, \\ [\gamma_2 \widetilde{\delta}_k, \ \widetilde{\delta}_k], & \eta_1 \leq \rho^k < \eta_2, \\ [\widetilde{\delta}_k, \ \infty], & \rho^k \geq \eta_2 \end{cases}$$

   and for   $x \in X^{k-1}$, set   $\delta_k(x) = \delta_{k-1}(x)$.

**Step 7: (Updating non-dominated solution set)**

   For   $x \in X^{k-1} \bigcup Y^k \bigcup \{x^{k+}\}$, compare the values   $F(x) = (f_1(x), f_2(x), \ldots, f_p(x))$

   and determine the non-dominated solution set   $X^k$.

**Step 8: (Iteration evaluation)**

   If   $\rho^k \geq \eta_1$   and   $X^k \setminus X^{k-1} \neq \varnothing$, then set   $k = k+1$   and go to Step 1.

   Else, set   $\delta_k(x^k) \in [\gamma_1 \widetilde{\delta}_k, \ \gamma_2 \widetilde{\delta}_k]$,   $x^{k+1} = x^k$,   $\widetilde{\delta}_{k+1} = \widetilde{\delta}_k$,   $k = k+1$, and go to Step 2.

## 4.  Convergence analysis

In this section we prove that the solutions produced by D-MOTR converge to Pareto critical points of the problem (2.1) on some assumptions. All the assumptions on the objective functions and the model functions are the same as the assumptions given in (Thomann, J; Eichfelder, G;, 2019). Proof strategies and steps for convergence analysis are performed similarly as in [6]. We assume that all objective functions   $f_i$,   $i \in \{1, \ldots, p\}$, are twice continuously differentiable and are bounded. We also assume that model functions   $m_i^k$   are quadratic and twice continuously differentiable for every index   $i \in \{1, \ldots, p\}$   and for every iteration   $k \in \mathbf{N}$. For every model, we assume that   $m^k(x^k) = F(x^k)$   holds for all   $k \in \mathbf{N}$.

**Assumption 1.** For every index $i \in \{1, ..., p\}$, the Hessian of the objective function $f_i$ is bounded, that is, there exists a constant $\kappa_{hf_i} > 1$ such that $||\nabla^2 f_i(x)|| \leq \kappa_{hf_i} - 1$ for all $x \in \mathbf{R}^n$.

**Assumption 2.** For every index $i \in \{1, ..., p\}$, the Hessian of the model function $m_i^k$ is uniformly bounded for all $k \in \mathbf{N}$, that is, there exists a constant $\kappa_{hm_i} > 1$ such that $||\nabla^2 m_i^k(x)|| \leq \kappa_{hm_i} - 1$ for all $x \in B_k$.

**Assumption 3.** For every index $i \in \{1, ..., p\}$, there exists a constant $\kappa_{e_i} > 0$ such that for all $k \in \mathbf{N}$ and $x \in B^k$

$$| f_i(x) - m_i^k(x) | \leq \kappa_{e_i} \widetilde{\delta}_k^{\,2}.$$

**Lemma 4.1.** ( (Thomann, J; Eichfelder, G;, 2019)). Suppose that
Assumption **1**,
Assumption **2**, and
Assumption **3** hold. Then there exists a constant $\kappa_{eg} > 0$ such that it holds

$$\left| \nabla f_i(x^k) - \nabla m_i^k(x^k) \right| \leq \kappa_{eg} \widetilde{\delta}_k$$

for all $i \in \{1, ..., p\}$ and $k \in \mathbf{N}$.

This lemma guarantees that the gradient of the model is very close to the gradient of the original function when the trust region radius is sufficiently small.

**Assumption 4.** There exists a constant $\kappa_\omega > 0$ such that for every $k \in \mathbf{N}$,

$$| \omega_m^k(x^k) - \omega(x^k) | \leq \kappa_\omega \widetilde{\delta}_k$$

holds.

**Theorem 4.2.** Suppose that Assumption 4 holds and $\delta_{tol} = 0$. If $\omega(x^k) \neq 0$ then the step 2 of algorithm terminates in finite time.

**Proof**. We suppose the step 2 of k-th iteration does not terminate in finite time. From the assumption 4, $| \omega_m^k(x^k) - \omega(x^k) | \leq \kappa_\omega \widetilde{\delta}_k$ holds. So it holds

$$\omega(x^k) \leq \omega_m^k(x^k) + \kappa_\omega \widetilde{\delta}_k \leq (1 + \kappa_\omega) \widetilde{\delta}_k.$$

The step 2 of k-th iteration repeats infinitely, so the right side of the inequality converges to zero. So $\omega(x^k) = 0$ holds and contrasts to assumptions of this lemma. $\qquad\square$

**Assumption 5.** ( (Thomann, J; Eichfelder, G;, 2019)). There exists a constant $\kappa_r \in (0, 1]$ such that for every $k \in \mathbf{N}$ with $x^k$ not Pareto critical for (3.1),

$$\frac{\min\limits_{i=1,\ldots,p} r_i^k}{\max\limits_{i=1,\ldots,p} r_i^k} \ge \kappa_r$$

holds.

**Lemma 4.3.** ( (Thomann, J; Eichfelder, G;, 2019)). Let $k \in \mathbf{N}$ be an iteration and $x^k$ be not a Pareto critical for (3.1). Let $g^k : \mathbf{R}^n \to \mathbf{R}$ be the quadratic function which is $g^k(x) := \sum_{i=1}^p \alpha_i m_i^k(x)$ with $\alpha_i \ge 0$ , $i \in \{1, \ldots, p\}$ . We suppose $g^k(x_c) = g^k(x^k + \bar{\tau}d_\omega) = \min\limits_{|t| \le \tilde{\delta}_k} g^k(x^k + td_\omega)$ with $d_\omega = -\nabla g^k(x^k) / ||\nabla g^k(x^k)||$ , and set $\beta_g^k := 1 + ||\nabla^2 g^k(x^k)||$ . Then it holds

$$g^k(x^k) - g^k(x_c) \ge \frac{1}{2} ||\nabla g^k(x^k)|| \min \left\{ \frac{||\nabla g^k(x^k)||}{\beta_g^k}, \tilde{\delta}_k \right\}.$$

**Lemma 4.4.** Suppose that
Assumption **2**,
Assumption **4** and
Assumption **5** hold. Let $x^{k+}$ be the solution of (3.2) and $\beta_m^k = \max\limits_{i=1,\ldots,p} ||\nabla^2 m_i^k(x^k)|| + 1$ .
Then for each $k \in \mathbf{N}$ , there exists an index $j_k \in \mathbf{N}$ such that

$$m_i^k(x^k) - m_i^k(x^{k+}) \ge \kappa_r \left( \frac{1}{2} \right)^{j_k} \omega_m^k(x^k) \min \left\{ \frac{\omega_m^k(x^k)}{\beta_m^k}, \tilde{\delta}_k \right\}.$$

**Proof**. Let $\left( t^{k+}, x^{k+} \right) \in \mathbf{R}^{1+n}$ be the solution of (3.2). First, let $x^k$ be not a Pareto critical for (3.1). Then from Lemma 3.3 and
Remark **3.1**, $t^{k+} \in [-1, 0)$ and $r_i^k = m_i^k(x^k) - \min\limits_{x \in B^k} m_i^k(x) > 0$ for all $i \in \{1, \ldots, p\}$ .
According to the constraints of (3.2),

$$m_i^k(x^k) - m_i^k(x^{k+}) \ge -t^{k+}r_i^k > 0$$

for all $i \in \{1, \ldots, p\}$ . Then

$$-t^{k+} = \left| t^{k+} \right| \le \frac{m_i^k(x^k) - m_i^k(x^{k+})}{r_i^k}. \tag{4.1}$$

Now, let $d_\omega \in \arg\min_{\|d\| \le 1} \max_{i=1,\ldots,p} \nabla m_i^k(x^k)^T d$ . Then from
Lemma **2.7**, there exists constants $\alpha_i \in [0, 1]$ , $i \in \{1, \ldots, p\}$ , with $\sum_{i=1}^p \alpha_i = 1$ and $\mu \geq 0$ such that $\|d_\omega\| = 1$. Let $x_c = x^k + \tau d_\omega$ with $|\tau| \leq \tilde{\delta}_k$ be the corresponding point from

Lemma **4.3** for the function $g^k(x) = \sum_{i=1}^p \alpha_i m_i^k(x)$.

Furthermore,

$$\beta_g^k = 1 + ||\nabla^2 g^k(x^k)|| \leq 1 + \sum_{i=1}^p \alpha_i ||\nabla^2 m_i^k(x^k)|| \leq 1 + \max_{i=1,\,...,\,p} ||\nabla^2 m_i^k(x^k)|| \models \beta_m^k$$

and according to
Lemma **4.3**,

$$g^k(x^k) - g^k(x_c) \geq \frac{1}{2} ||\nabla g^k(x^k)|| \min\left\{ \frac{||\nabla g^k(x^k)||}{\beta_m^k}, \tilde{\delta}_k \right\}.$$

$$(4.2)$$

Since $x_c \in B^k$ and $d_\omega$ is a descent direction for (3.1) and

Lemma **2.7**, there exists a constant $t$ such that $(t, x_c)$ is feasible for (3.2). Then there exists a smallest constant $t_c$ such that $(t_c, x_c)$ is feasible for (3.2) and

$$-t_c = |t_c| \models \min_{i=1,\,...,\,p} \frac{m_i^k(x^k) - m_i^k(x_c)}{r_i^k} \geq \frac{\min\limits_{i=1,\,...,\,p} \left( m_i^k(x^k) - m_i^k(x_c) \right)}{\max\limits_{i=1,\,...,\,p} r_i^k}.$$

As $t^{k+}$ is the minimal value of (3.2), $|t_c| \leq |t^{k+}|$.
From (4.1),

$$m_i^k(x^k) - m_i^k(x^{k+}) \geq \min_{i=1,\,...,\,p} \left( m_i^k(x^k) - m_i^k(x_c) \right) \cdot \kappa_r. \qquad (4.3)$$

Since $\sum_{i=1}^p \alpha_i = 1$ and $(t_c, x_c)$ is feasible for (3.2),

$$g^k(x^k) - g^k(x_c) = \sum_{i=1}^p \alpha_i \left[ m_i^k(x^k) - m_i^k(x_c) \right] \geq \min_{i=1,\,...,\,p} \left( m_i^k(x^k) - m_i^k(x_c) \right) > 0.$$

Then from (4.2), there exists an index $j_k \in \mathbf{N}$ such that

$$\min_{i=1,\,...,\,p} \left( m_i^k(x^k) - m_i^k(x_c) \right) \geq \left( \frac{1}{2} \right)^{j_k} ||\nabla g^k(x^k)|| \min\left\{ \frac{||\nabla g^k(x^k)||}{\beta_m^k}, \tilde{\delta}_k \right\}$$

. Therefore, from (4.3),

$$m_i^k(x^k) - m_i^k(x^{k+}) \geq \left( \frac{1}{2} \right)^{j_k} \left\| \sum_{i=1}^p \alpha_i \nabla m_i^k(x^k) \right\| \min\left\{ \frac{\left\| \sum_{i=1}^p \alpha_i \nabla m_i^k(x^k) \right\|}{\beta_m^k}, \tilde{\delta}_k \right\} \cdot \kappa_r \qquad (4.4)$$

for every $k \in \mathbf{N}$ with $x^k$ being not a Pareto critical point.

If $x^k$ is a Pareto critical point for (3.1), then $\omega_m^k(x^k) = 0$ and $d_\omega = 0$. Therefore,

$\sum_{i=1}^{p} \alpha_i \nabla m_i^k (x^k) = 0$. As $x^{k+}$ is the solution of (3.2), $m_i^k(x^k) - m_i^k(x^{k+}) \geq 0$ and (4.4) is also satisfied.

According to

Lemma **2.7**, $\omega_m^k (x^k) \leq \left\| \sum_{i=1}^{p} \alpha_i \nabla m_i^k (x^k) \right\|$. So,

$$m_i^k (x^k) - m_i^k (x^{k+}) \geq \left( \frac{1}{2} \right)^{j_k} \left\| \sum_{i=1}^{p} \alpha_i \nabla m_i^k (x^k) \right\| \min \left\{ \frac{\left\| \sum_{i=1}^{p} \alpha_i \nabla m_i^k (x^k) \right\|}{\beta_m^k}, \tilde{\delta}_k \right\} \cdot \kappa_r$$

$$\geq \left( \frac{1}{2} \right)^{j_k} \kappa_r \omega_m^k (x^k) \min \left\{ \frac{\omega_m^k (x^k)}{\beta_m^k}, \tilde{\delta}_k \right\}$$

for every iteration $k \in \mathbf{N}$.    $\square$

**Assumption 6.** There exists a constant $\kappa_m \in (0, 1)$ such that for all $k \in \mathbf{N}$,

$$m_i^k (x^k) - m_i^k (x^{k+}) \geq \kappa_m \omega(x^k) \min \left\{ \frac{\omega(x^k)}{\beta_m^k}, \tilde{\delta}_k \right\}$$

with $\beta_m^k = \max_{i=1, \ldots, p} \left\| \nabla^2 m_i^k (x^k) \right\| + 1$.

**Remark 4.5.** From the Assumption **2**,

$$\beta_m^k = \max_{i=1, \ldots, p} || \nabla^2 m_i^k (x^k) || + 1 \leq \max_{i=1, \ldots, p} \kappa_{hm_i}$$

for all $k \in \mathbf{N}$.

Now consider the following sets of iterations.

$$S := \left\{ k \in \mathbf{N} \,\middle|\, \rho^k = \min_{i=1, \ldots, p} \frac{f_i(x^k) - f_i(x^{k+})}{m_i^k(x^k) - m_i^k(x^{k+})} \geq \eta_1 \right\}; \text{ set of all successful iterations}$$

$$V := \left\{ k \in \mathbf{N} \,\middle|\, \rho^k = \min_{i=1, \ldots, p} \frac{f_i(x^k) - f_i(x^{k+})}{m_i^k(x^k) - m_i^k(x^{k+})} \geq \eta_2 \right\}; \text{ set of all very successful iterations}$$

$$U := \left\{ k \in \mathbf{N} \,\middle|\, \rho^k = \min_{i=1, \ldots, p} \frac{f_i(x^k) - f_i(x^{k+})}{m_i^k(x^k) - m_i^k(x^{k+})} < \eta_1 \right\}; \text{ set of all unsuccessful iterations}$$

**Lemma 4.6.** Suppose that Assumption **1**, Assumption **2**,

Assumption **3**,
Assumption **4**,
Assumption **5** and
Assumption **6** hold, $x^k$ is not a Pareto critical for (2.1) and

$$\tilde{\delta}_k \leq \frac{\kappa_m (1 - \eta_2) \omega_m^k(x^k)}{\kappa_e} . \qquad (4.5)$$

Here, $\kappa_e := \max\limits_{i=1,\,\ldots,\,p} \max\{\kappa_{e_i}, \kappa_{hm_i}\} > 0$. Then $k \in V$, that is iteration $k$ is very successful.

**Proof**. As $x^k$ is not a Pareto critical for (2.1), according to
Lemma **2.6** it holds $\omega_m^k(x^k) > 0$. By (4.5),

$$\tilde{\delta}_k \leq \frac{\kappa_m (1 - \eta_2) \omega_m^k(x^k)}{\kappa_e} .$$

Due to $\eta_2$, $\kappa_m \in (0, 1)$ and $\kappa_m (1 - \eta_2) < 1$,

$$\frac{\kappa_m (1 - \eta_2) \omega_m^k(x^k)}{\kappa_e} < \frac{\omega_m^k(x^k)}{\kappa_e} .$$

By the definition $\kappa_e$, $\kappa_e \geq \max\limits_{i=1,\,\ldots,\,p} \kappa_{hm_i}$ and

$$\frac{\omega_m^k(x^k)}{\kappa_e} \leq \frac{\omega_m^k(x^k)}{\max\limits_{i=1,\,\ldots,\,p} \kappa_{hm_i}} .$$

According to
Remark **4.5**,

$$\frac{\omega_m^k(x^k)}{\max\limits_{i=1,\,\ldots,\,p} \kappa_{hm_i}} \leq \frac{\omega_m^k(x^k)}{\beta_m^k} .$$

So it follows

$$\tilde{\delta}_k \leq \frac{\omega_m^k(x^k)}{\beta_m^k} .$$

According to
Assumption **6**, it holds

$$m_i^k(x^k) - m_i^k(x^{k+}) \geq \kappa_m \omega_m^k(x^k) \min\left\{\frac{\omega_m^k(x^k)}{\beta_m^k}, \tilde{\delta}_k\right\} \geq \kappa_m \omega_m^k(x^k) \tilde{\delta}_k . \qquad (4.6)$$

From the interpolation condition $m_i^k(x^k) = f_i^k(x^k)$, it holds

$$\left| \rho_i^k - 1 \right| = \left| \frac{f_i(x^k) - f_i(x^{k+})}{m_i^k(x^k) - m_i^k(x^{k+})} - 1 \right| = \left| \frac{m_i^k(x^k) - f_i(x^{k+})}{m_i^k(x^k) - m_i^k(x^{k+})} - 1 \right|$$

$$= \left| \frac{m_i^k(x^{k+}) - f_i^k(x^{k+})}{m_i^k(x^k) - m_i^k(x^{k+})} \right| .$$

From
Assumption **3**, (4.6) it follows

$$\left| \frac{m_i^k(x^{k+}) - f_i^k(x^{k+})}{m_i^k(x^k) - m_i^k(x^{k+})} \right| \le \frac{\kappa_{e_i} \widetilde{\delta}_k{}^2}{\kappa_m \omega_m^k(x^k) \widetilde{\delta}_k} \le \frac{\max\limits_{i=1,\dots,p} \kappa_{e_i} \widetilde{\delta}_k}{\kappa_m \omega_m^k(x^k)} \le \frac{\kappa_e \widetilde{\delta}_k}{\kappa_m \omega_m^k(x^k)} \le 1 - \eta_2 \,.$$

This implies $\rho_i^k \ge \eta_2$ and therefore $\rho^k \ge \eta_2$ and $k \in V$. $\qquad\square$

**Definition 4.7.** ( (Liuzzi, G; Lucidi, S; Rinaldi, F;, 2016)) We define a linked sequence as a sequence $\{x^{k_i}\}$ such that for any $i = 1, 2, \dots,$ $x^{k_{i+1}}$ is generated by $x^{k_i}$ in D-MOTR.

**Lemma 4.8.** Suppose that
Assumption **1**,
Assumption **2**,
Assumption **3**,
Assumption **4**,
Assumption **5** and
Assumption **6** hold. Let $\{x^{k_i}\}$ be any linked sequence of algorithm and $K$ ($k_i \in K \subset \mathbf{N}$) be set of iterations. If $\lim\limits_{i \to +\infty} \widetilde{\delta}_{k_i} = 0$, then $\liminf\limits_{i \to +\infty} \omega_m^{k_i}(x^{k_i}) = 0$.

**Proof**. Assume $\liminf\limits_{i \to +\infty} \omega_m^{k_i}(x^{k_i}) = \kappa_{b\omega} > 0$. Then there exists a $k_i \in K$ for any $\varepsilon > 0$ such that $\widetilde{\delta}_{k_i} < \varepsilon$. Consider

$$\varepsilon = \frac{\gamma_1 \kappa_m \kappa_{b\omega}(1 - \eta_2)}{\kappa_e} \,.$$

Let $k_s$ be the first iteration such that $k_{s-1} \in K'$ and $\widetilde{\delta}_{k_s} < \varepsilon$. Then

$$\widetilde{\delta}_{k_{s-1}} \ge \varepsilon > \widetilde{\delta}_{k_s} \,. \tag{4.7}$$

According to the step 6 of D-MOTR, it holds $\gamma_1 \widetilde{\delta}_{k_{s-1}} \le \widetilde{\delta}_{k_s}$. So,

$$\widetilde{\delta}_{k_{s-1}} \le \frac{\widetilde{\delta}_{k_s}}{\gamma_1} < \frac{\varepsilon}{\gamma_1} = \frac{\kappa_m \kappa_{b\omega}(1 - \eta_2)}{\kappa_e} \le \frac{\kappa_m \omega_m^{k_{s-1}}(x^{k_{s-1}})(1 - \eta_2)}{\kappa_e} \,.$$

Because of $\omega_m^{k_{s-1}}(x^{k_{s-1}}) \ge \kappa_{b\omega} > 0$ and
Lemma **2.6**, $x^{k_{s-1}}$ is not a Pareto critical for (2.1). Therefore, from
Lemma **4.6**, $k_{s-1} \in V$. So according to the step of D-MOTR, $\widetilde{\delta}_{k_{s-1}} \le \widetilde{\delta}_{k_s}$. This contradicts equation (4.7). $\qquad\square$

**Lemma 4.9.** Suppose that
Assumption **1**,
Assumption **2**,

Assumption **3**,

Assumption **4**,

Assumption **5** and

Assumption **6** hold and $\{x^{k_i}\}$ ( $k_i \in K \subset \mathbf{N}$ ) be any linked sequence of D-MOTR. If the set of successful iterations, $K \cap S$ is finite then there exists $k_s \in K$ with $x^{k_{s+1}}$ be a Pareto critical for (2.1) and $x^{k_i} = x^{k_{s+1}}$ for any $k_i \geq k_{s+1}$ ( $k_i \in K$ ).

**Proof**. Let $k_s$ be the last successful iteration in $K$. Then it holds $\rho^{k_i}(x) < \eta_1$ for every $k_i > k_s$. Step 8 of the D-MOTR ensures $x^{k_{i+1}} = x^{k_i}$ and $\delta_{k_{i+1}}(x^{k_{i+1}}) \leq \gamma_2 \delta_{k_i}(x^{k_i})$. So

$$\lim_{i \to +\infty} \delta_{k_i}(x^{k_i}) = \lim_{i \to +\infty} \delta_{k_i}(x^{k_{s+1}}) = 0 . \tag{4.8}$$

Assume that $x^{k_{s+1}}$ is not a Pareto critical for (2.1). From (4.8), there exists $\hat{k} \geq k_{s+1}$ ( $\hat{k} \in K$ ) such that

$$\widetilde{\delta}_{\hat{k}} = \delta_{\hat{k}}(x^{\hat{k}}) \leq \frac{\kappa_m (1-\eta_2)\omega_m^{k_{s+1}}(x^{k_{s+1}})}{\kappa_e} = \frac{\kappa_m (1-\eta_2)\omega_m^{\hat{k}}(x^{\hat{k}})}{\kappa_e} .$$

Then

Lemma **4.6** implies that every $\hat{k}$ is a very successful iteration. This contradict $k_s$ is the last successful iteration in $K$. Hence $x^{k_{s+1}}$ is a Pareto critical for (2.1) and $\omega_m^{k_{s+1}}(x^{k_{s+1}}) = 0$. And from Assumption 4 and (4.8), $\omega(x^k) \leq \omega_m^k(x^k) + \kappa_\omega \widetilde{\delta}_k$ and $\lim_{i \to +\infty} \omega(x^{k_i}) = 0$.                                       $\square$

**Lemma 4.10.** Suppose that

Assumption **1**,

Assumption **2**,

Assumption **3**,

Assumption **4**,

Assumption **5** and

Assumption **6** hold and $\{x^{k_i}\}$ ( $k_i \in K \subset \mathbf{N}$ ) be any linked sequence of D-MOTR. If the set of successful iteration $K \cap S$ is infinite, then

$$\liminf_{i \to +\infty} \omega(x^{k_i}) = 0 .$$

**Proof**. From

Remark **4.5** it follows

$$\beta_m^{k_i} \leq \max_{j=1,\ldots,p} \kappa_{hm_j} \leq \max_{j=1,\ldots,p} \max\{\kappa_{hm_j}, \kappa_{e_i}\} = \kappa_e$$

for every iteration $k_i \in K^r$.

Since for all $k_i \in K^r$,

$$\min_{j=1, \dots, p} \frac{f_j(x^{k_i}) - f_j(x^{k_i+})}{m_j^{k_i}(x^{k_i}) - m_j^{k_i}(x^{k_i+})} \geq \eta_1,$$

it follows $f_j(x^{k_i}) - f_j(x^{k_i+}) \geq \eta_1(m_j^{k_i}(x^{k_i}) - m_j^{k_i}(x^{k_i+}))$.

Now, let consider the hyper-rectangle $[F(x); R_{HV}]$ which has $R_{HV} = [\max f_1(x), \max f_2(x), \dots, \max f_p(x)]$ and $F(x), x \in X^{k_i}$, as edge points of diagonal. Then, let $HV^{k_i}$ be the volume of

$$\bigcup_{x_j \in X^{k_i}} [F(x_j); R_{HV}],$$

which is the union of all hyper-rectangles produced by non-dominated points of $X^{k_i}$ and $R_{HV}$.

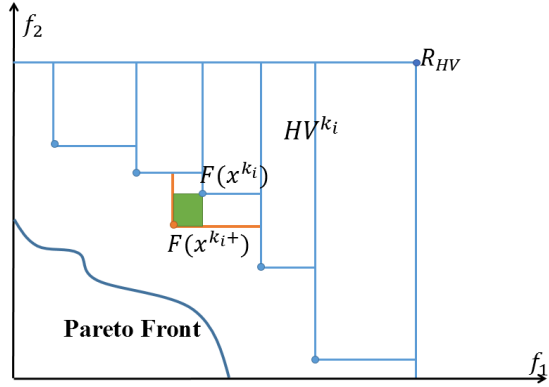

Fig 4. Pareto front and $HV^{k_i}$

From
Assumption **6**,

$$HV^{k_{i+1}} - HV^{k_i} \geq \prod_{j=1}^{p} \left( f_j(x^{k_i}) - f_j(x^{k_i+}) \right)$$

$$\geq \eta_1^{p} \cdot \prod_{j=1}^{p} \left( m_j^k(x^{k_i}) - m_j^k(x^{k_i+}) \right)$$

$$\geq \left( \eta_1 \kappa_m \omega_m^{k_i}(x^{k_i}) \min \left\{ \frac{\omega_m^{k_i}(x^{k_i})}{\beta_m^{k_i}}, \widetilde{\delta}_{k_i} \right\} \right)^{p}$$

$$\geq \left( \eta_1 \kappa_m \varepsilon \min \left\{ \frac{\varepsilon}{\kappa_e}, \widetilde{\delta}_{k_i} \right\} \right)^{p} \tag{4.9}$$

Since $S$ is infinite, $\{HV^k\}$ is monotonically increasing and bounded and it holds

$HV^{k_{i+1}} - HV^{k_i} \geq 0$ , (4.9) converges to zero. Hence, $\lim\limits_{i \to +\infty} \widetilde{\delta}_{k_i} = 0$ . So, from

Lemma **4.8** and

Assumption **4**,

$$\liminf\limits_{i \to +\infty} \omega(x^{k_i}) = 0 .$$ □

**Theorem 4.11.** Suppose that

Assumption **1**,

Assumption **2**,

Assumption **3**,

Assumption **4**,

Assumption **5** and

Assumption **6** hold. Then,

$$\liminf\limits_{i \to +\infty} \omega(x^{k_i}) = 0$$

for any linked sequence $\{x^{k_i}\}$ of D-MOTR.

## 5. Numerical evaluations.

### 5.1. *Performance measures.*

In this paper, we evaluate problems with box-constraints. We handle the constraints by reducing trust region of each points to be included in feasible region. In multiobjective optimization, it is important how close the generated points are to Pareto optimality and how well the images of set of these solutions are distributed over the Pareto front.

In multi-objective optimization, there are several indicators that evaluate the accuracy of the solutions, such as GD, GDp, ONVG, RNI, and so on, but we use GD to be inconvenient, since we experiment only with at most three objective functions.

Let $N = \{x_1, x_2, \ldots, x_v\}$ denote the set of produced solutions by an algorithm.

First, we see the GD which is the closeness measure of produced solutions. GD is a well-known measure evaluating the accuracy of set of approximated solutions by calculating the average distance between the Pareto front and produced solutions (VeldhuizenD, 1999). Let $F(\hat{x}_j)$, $j = 1, \ldots, M$ denote the Pareto efficient point which is closest to $F(y_j)$. Then, GD is computed as

$$GD = \frac{\sqrt{\sum\limits_{j=1}^{M} \left\| F(y_j) - F(\hat{x}_j) \right\|^2}}{M} .$$

So, if the differences between generated solutions and Pareto efficiency are smaller, the GD goes to zero.

There are several criteria measures such as HV, IGD, Minimal Spacing, Entropy

Measure and so on for diversity of solutions (Farhang-Mehr, A; Azarm, S;, 2002; Wang, H; Jin, Y; Yao, X;, 2017; Li, M; Yang, S; Liu, X;, 2014). Among them, Hypervolume is one of the best measures to evaluate the diversity in multiobjective optimization. And a new algorithm is produced to calculate the Hypervolume with higher speed (Lacour, R; Klamroth, K; Fonseca, C M;, 2016). Hypervolume is equal to the volume enclosed by the generated points and a reference point, say, $R_{HV} = (r_1, r_2, \cdots, r_p)$, which must be dominated by all produced points. Often, we set $R_{HV} = (r_1, r_2, \cdots, r_p)$ as $R_{HV} = (\max f_1(x), \max f_2(x), \cdots, \max f_p(x))$. Let $[F(x); R_{HV}]$ which is mentioned in Lemma 4.9, denote the hyper-rectangle. Then, Hypervolume is calculated by

$$\bigcup_{x_j \in X^{(k)}} [F(x_j); R_{HV}].$$

## 5.2.   *Illustrations of the algorithm*

### 5.2.1.  *Comparison with BIMADS*

We test the D-MOTR algorithm with the following problem.

$$\min \left( \begin{array}{l} f_1(x_1, \ldots, x_4) = 1 - \exp\left( \sum_{i=1}^{4} (x_i - 1/2)^2 \right) \\ f_2(x_1, \ldots, x_4) = 1 - \exp\left( -\sum_{i=1}^{4} (x_i - 1/2)^2 \right) \end{array} \right)$$

Left part is results for BIMADS and right one is for D-MOTR in Fig 5. Here, the curve is the Pareto front and green points stand for the evaluated points before. Blue points are non-dominated points. To implement BIMADS, we used the MATLAB function patternsearch with MADS. For the comparison, we set the initial parameters of BIMADS as follows:

- InitialMeshSize = 0.5
- MeshExpansion = 2
- MeshContraction = 0.5
- TolMesh = 1e-7
- MaxFunEvals = 50

After some iterations of BIMADS algorithm, since the solutions of single objective optimization don't improve non-dominated points, the reference point is fixed. To avoid this, we assigned bigger weight to the selected reference point than the others in each iteration. For all our effort, the BIMADS algorithm does not give sufficient results, because we may not use MADS for single-objective optimization and initiated proper parameters.

Parameters for D-MOTR are set to $x^0 = (0.1, -0.1, 0.1, -0.1)$, $\delta_0(x^0) = 1$, $\delta_{tol} = 0.05$, $\eta = 0.5$, $\gamma_0 = 0.7$, $\gamma_1 = 0.5$, $\gamma_2 = 1$, $\sigma = 0.05$.

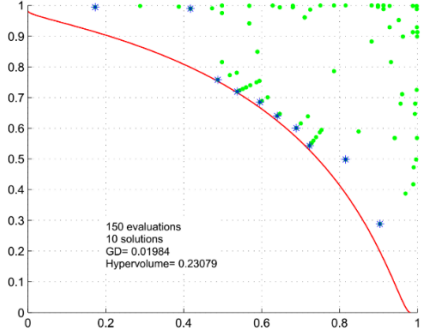

(a.1). BIMADS with 150 evaluations

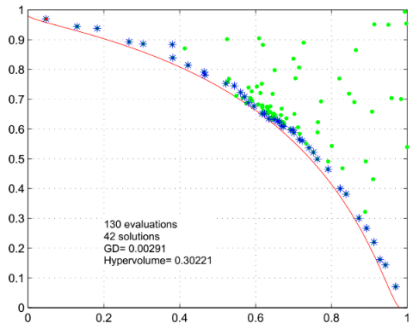

(b.1). D-MOTR with 130 evaluations

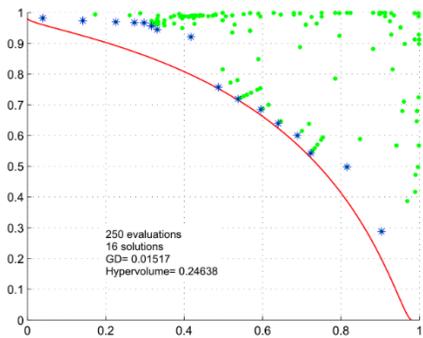

(a.2). BIMADS with 250 evaluations

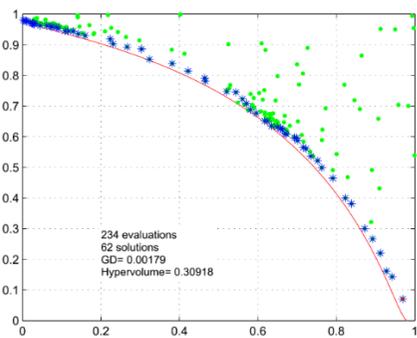

(b.2). D-MOTR with 234 evaluations

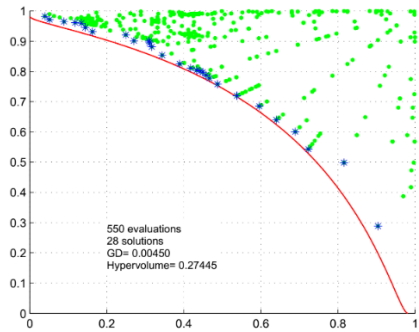

(a.3). BIMADS with 550 evaluations

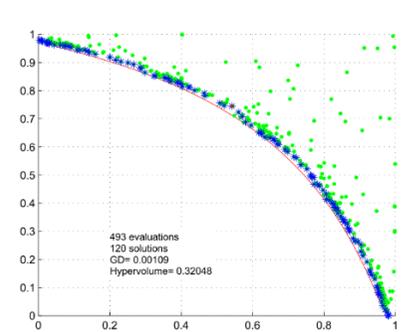

(b.3). D-MOTR with 493 evaluations

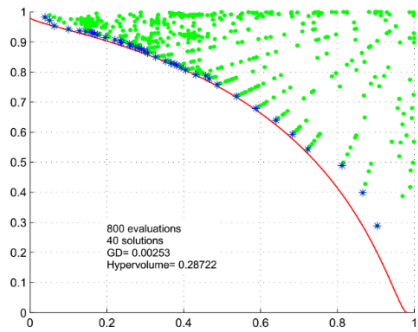
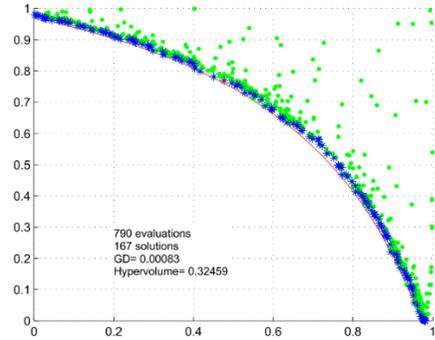

(a.4). BIMADS with 800 evaluations          (b.4). D-MOTR with 790 evaluations

Fig 5. D-MOTR compared to BIMADS

Left part in Fig 6 is progress of the D-MOTR algorithm in objective space. Here, the curve is the Pareto front and points around the curve are non-dominated points. The line at south-west in figure is projective line and points placed on it are projected points of non-dominated points. The other line normal to the projective line shows the direction of projection. Right part in figure shows the density function on projective line with the non-dominated points.

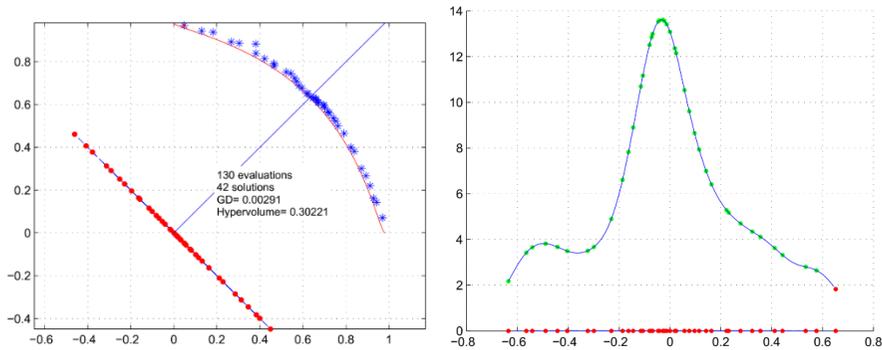

(a). After the first iteration.

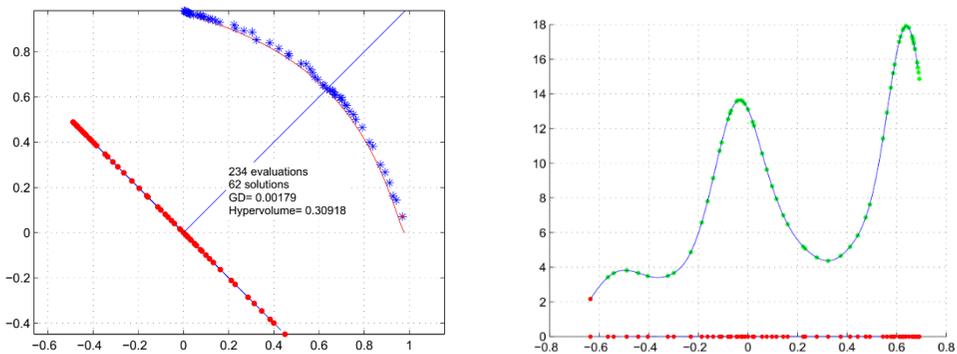

(b). After the second iteration.

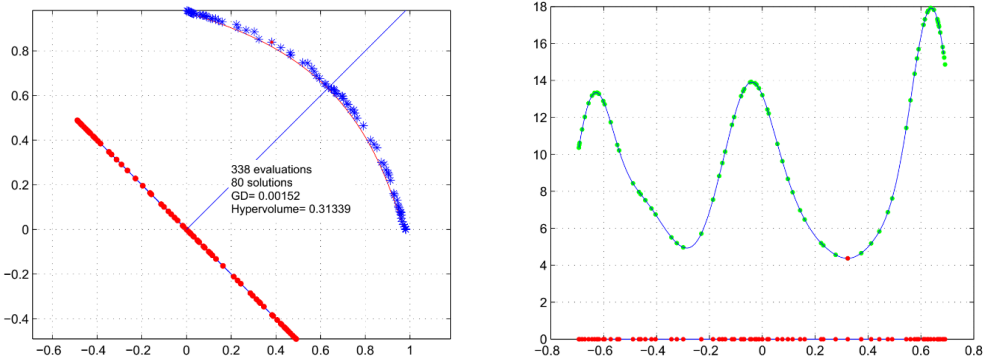

(c). After the third iteration.

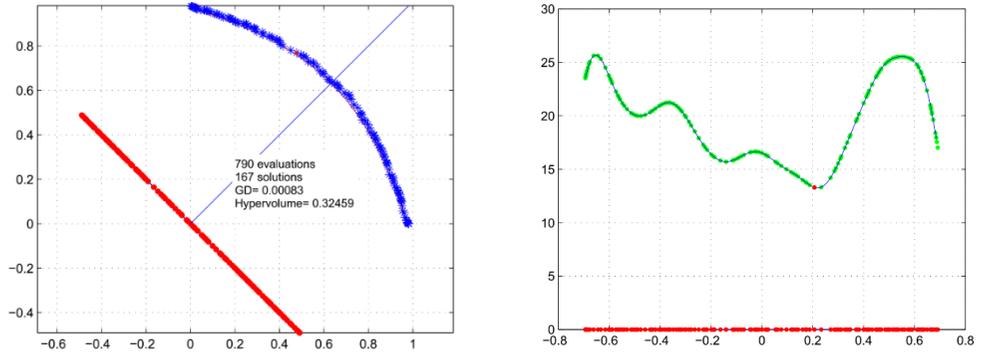

(d). After the ninth iteration.

Fig 6. Progress of the D-MOTR algorithm.

We can notice that the value of density function is small in sparse part and is big in dense part. The value which has the smallest density function among the non-dominated point is plotted with red dot, and the point would be selected for reference point in next iteration.

### 5.2.2. *Simple illustration with DTLZ2.*

We test the D-MOTR algorithm with a tri-objective optimization problem with three variables, DTLZ2, adopted from (Abraham, A; Jain, L; Goldberg, R;, 2005), defined as follows:

$$
\min \quad
\begin{pmatrix}
f_1(\mathbf{x}) = \big(1 + g(\mathbf{x})\big)\cos(x_1 \cdot \pi / 2)\cos(x_2 \cdot \pi / 2) \\
f_2(\mathbf{x}) = \big(1 + g(\mathbf{x})\big)\cos(x_1 \cdot \pi / 2)\sin(x_2 \cdot \pi / 2) \\
f_3(\mathbf{x}) = \big(1 + g(\mathbf{x})\big)\sin(x_1 \cdot \pi / 2)
\end{pmatrix}
$$

$$\text{s.t.} \qquad g(\mathbf{x}) = (x_3 - 0.5)^2$$
$$0 \le x_i \le 1, \quad i = 1, 2, 3$$

The parameters of algorithm are $x^0 = (0.5, 0.5, 0.5)$, $\delta_0(x^0) = 1$, $\delta_{tol} = 0.05$, $\gamma_0 = 0.7$, $\gamma_1 = 0.5$, $\gamma_2 = 5$, $\sigma = 0.05$. Fig 7 shows the results of D-MOTR for DTLZ2. The left part of the figure shows the progress of algorithm on objective space. The mesh part is the Pareto Front and dots are non-dominated solutions obtained in the algorithm. The bog point is the reference point and would be selected for next iteration. The right part of the figure shows the density function on the hyper projective plane. The dots are the projected solutions on to the projective plane and the density function is shown as mesh grid. The big point is the reference point as above.

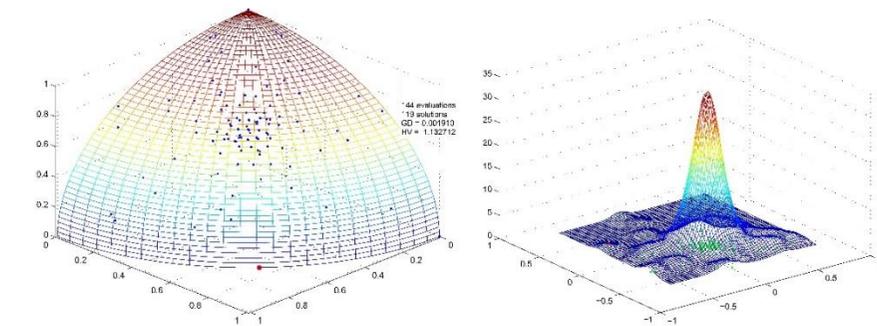

(a). First iteration of the algorithm

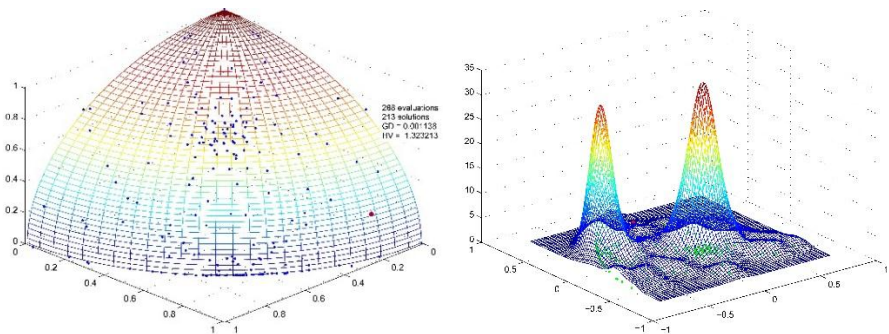

(b). Second iteration of the algorithm

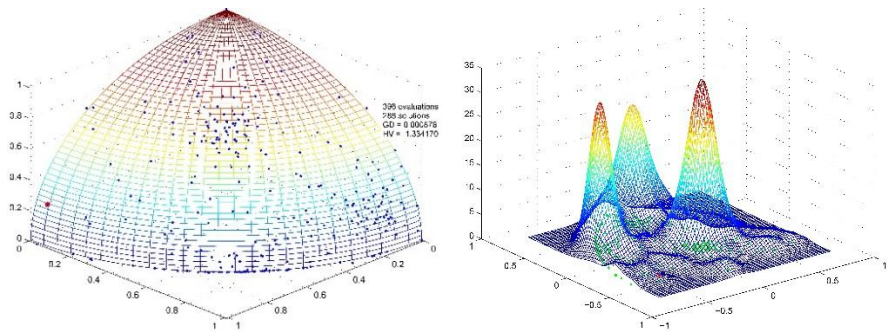

(c). Third iteration of the algorithm

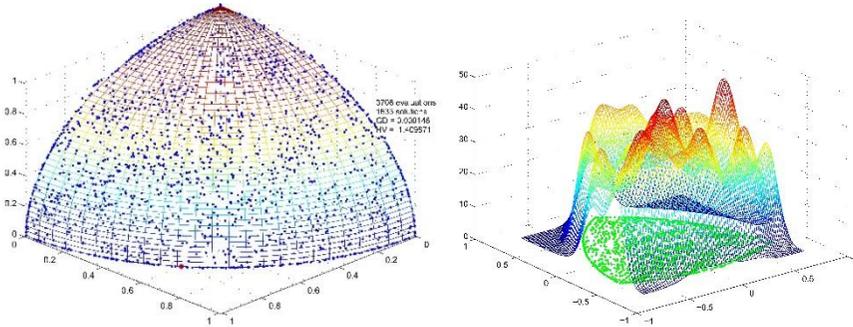

(d).40th iteration of the algorithm

Fig 7. Iterations of the D-MOTR with DTLZ2

In Fig 7.(d), the points are projected points of non-dominated points on to the hyper projective plane. We can find that the D-MOTR algorithm is generating Pareto approximated solutions effectively through this numerical evaluation.

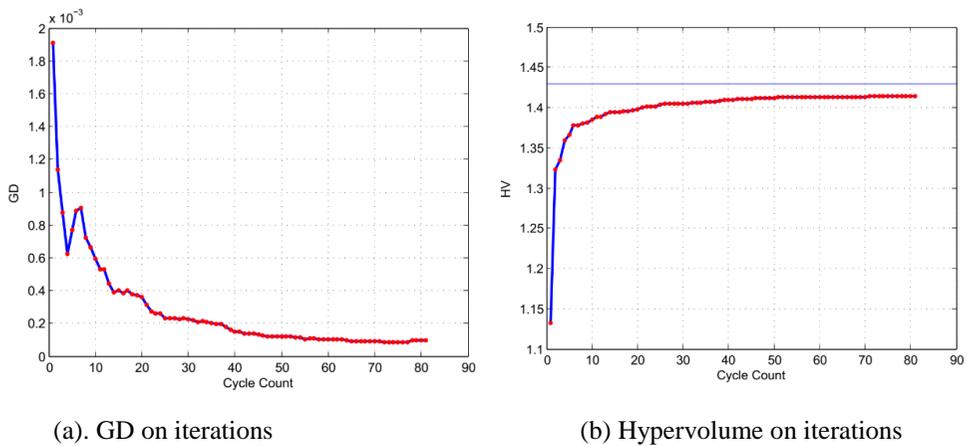

(a). GD on iterations               (b) Hypervolume on iterations

Fig 8. GD and Hypervolume on iterations of D-MOTR

Fig 8 shows the changes of GD and Hypervolume according to the iterations of algorithm. We can notice that GD approaches to zero while the iteration is increasing, in Fig (a). And Hypervolume approaches to a critical value, and this means that the set of non-dominated points is spanning to the Pareto front.

### 5.2.3. *Comet Problem*

We considered three-objective optimization problem from the literature (Abraham, A; Jain, L; Goldberg, R;, 2005), called Comet Problem.

$$\min \begin{pmatrix} f_1(x) = (1+x_3)(x_1^3 x_2^2 - 10x_1 - 4x_2) \\ f_2(x) = (1+x_3)(x_1^3 x_2^2 - 10x_1 + 4x_2) \\ f_3(x) = 3(1+x_3)x_1^2 \end{pmatrix}$$

$$\text{s.t.} \quad \begin{aligned} 1 &\le x_1 \le 3.5 \\ -2 &\le x_2 \le 2 \\ 0 &\le x_3 \le 1 \end{aligned}$$

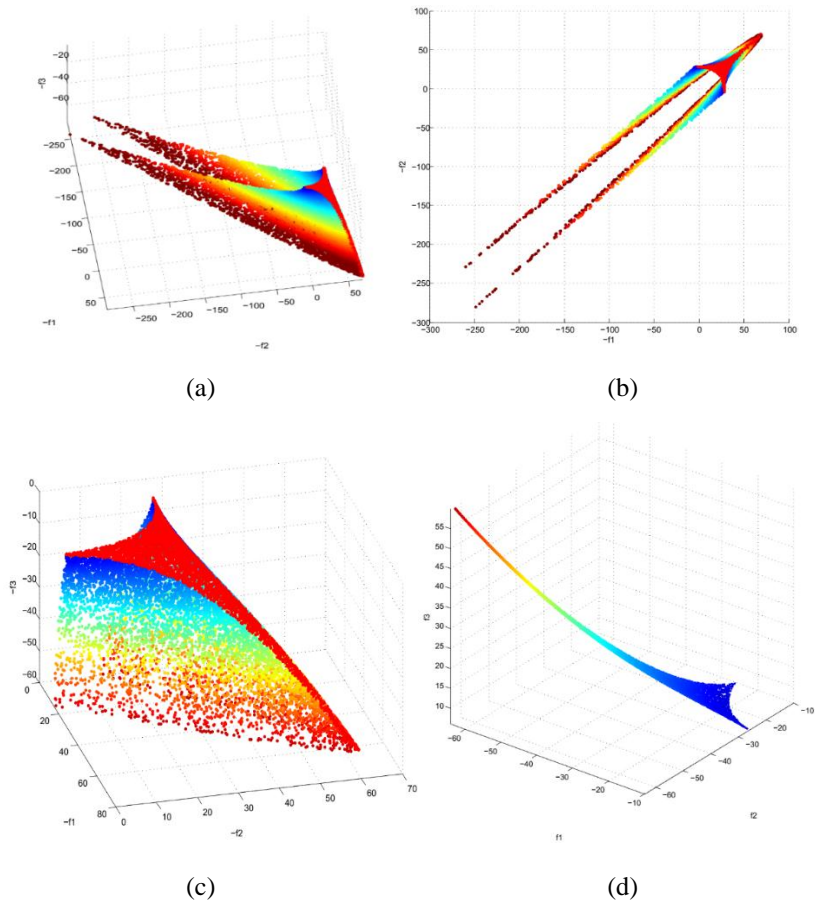

(a) (b)

(c) (d)

Fig 9. Area of Pareto efficiencies in Comet problem.

In this problem, Pareto efficiencies are concentrated on very small area compared to the feasible region. And the Pareto front is consisted of narrow part and wide part. So, we should visit not only the wide part of the Pareto front, but also the narrow part. Because of these properties, Comet problem is known as one of the test problem for performance of multiobjective optimization algorithm.

Fig 9 shows the area of Pareto efficiencies of Comet problem in objective space. The narrow part in Fig (a)-(c) is the area of Pareto efficiencies and Fig (d) shows only the area.

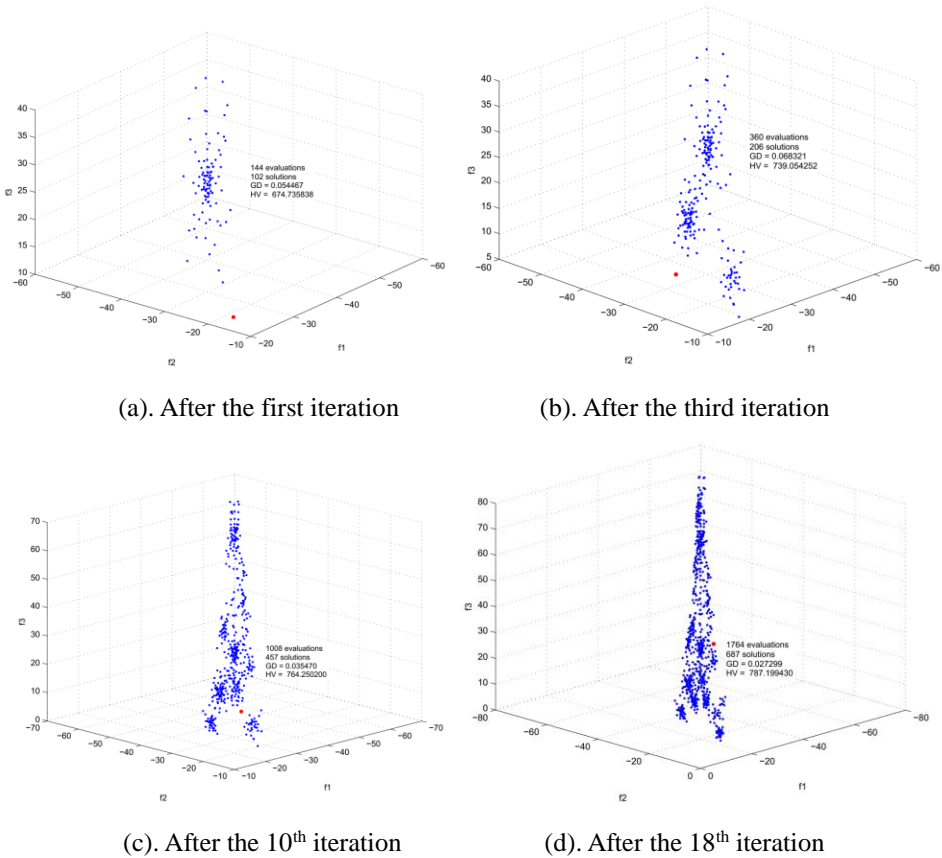

(a). After the first iteration             (b). After the third iteration

(c). After the 10th iteration             (d). After the 18th iteration

Fig 10. Process of performance with Comet problem through the 18th iteration.

We can notice that the solutions are spreading over the Pareto front while the iteration of algorithm is increasing. You can see that the algorithm visits the narrow and wide part of the Pareto front, coincidently.

### 5.2.4. *Disconnected Pareto front.*

D-MOTR algorithm is designed to search for a connected local Pareto front. However, we test the algorithm with DTLZ7 from (EichfelderG., 2008) in which the Pareto front is not connected.

$$\min \quad \begin{pmatrix} f_1 = x_1 \\ f_2 = x_2 \\ f_3 = [1 + g(x)] \left( 3 - \sum_{i=1}^{2} \left( \frac{x_i}{1 + g(x)} \left( 1 + \sin(3\pi x_i) \right) \right) \right) \end{pmatrix}$$

s.t. $\quad g(x) = 1 + \dfrac{9}{2} x_3, \ x_i \in [0, 1], \ i = 1, 2, 3,$

Pareto front of this problem is disconnected and consists of four connected parts.

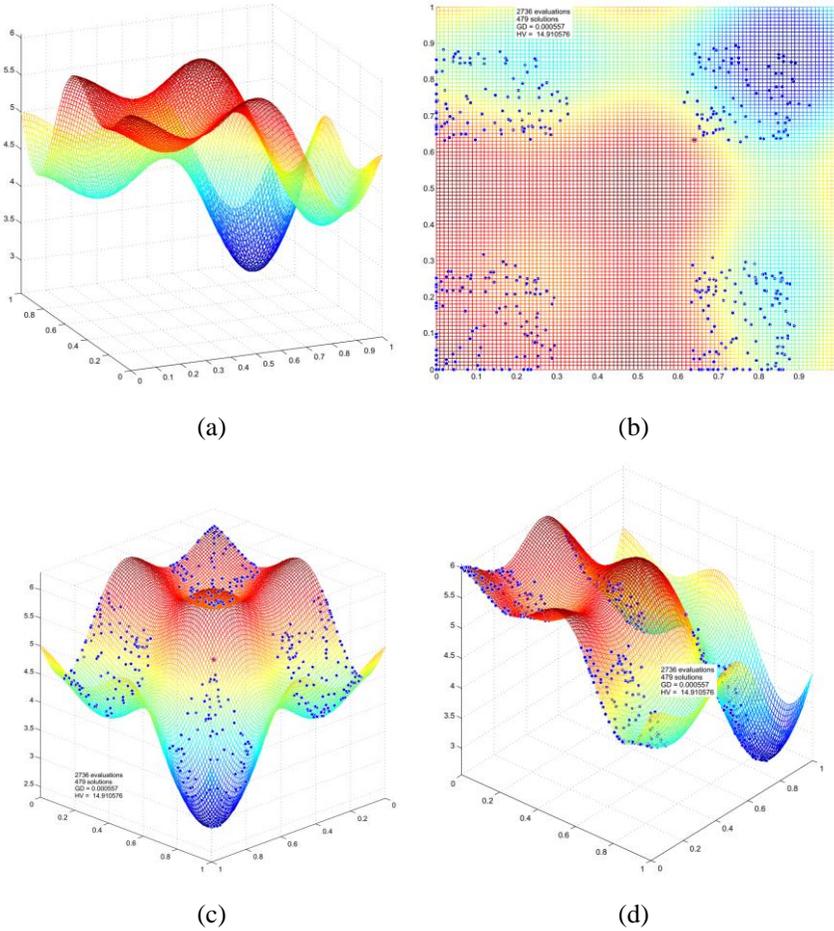

(a)                          (b)

(c)                          (d)

Fig 11. Performance of D-MOTR for DTLZ7 with a disconnected Pareto front.

And set of efficient solutions is a subset of the set

$$B := \left\{ y \in R^3 \mid y_1, \ y_2 \in [0, 1], \ y_3 = 2 \cdot \left( 3 - \sum_{i=1}^{2} \left( \frac{y_i}{2} (1 + \sin(3\pi y_i)) \right) \right) \right\},$$

which is plotted in Fig 11 (a).

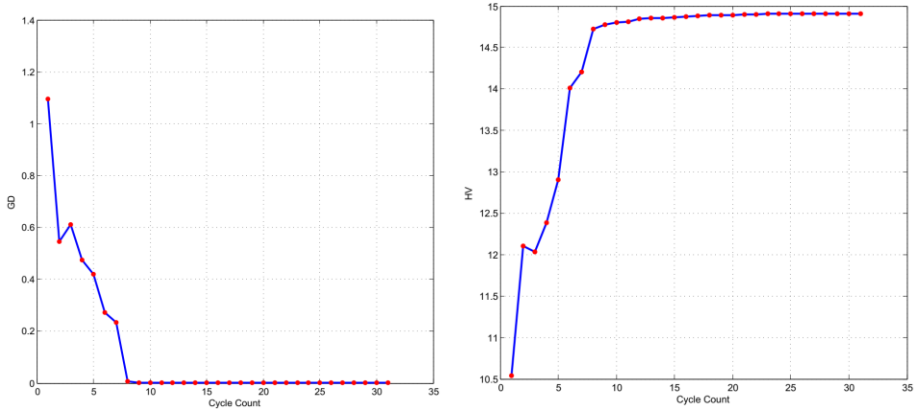

Fig 12. GD and Hypervolume with respect to the cycle count of D-MOTR

## 6. Concluding remarks

In this paper, we proposed a trust region algorithm for the Pareto front approximation in black-box multiobjective optimization. In the algorithm, we employed the density function to select a reference point in the set of non-dominated points and explored area around this point. This ensured uniformity of solutions distribution even for problems with more than two objective functions. We also proved that the iteration points of the algorithm converge to Pareto critical points. Finally, we presented numerical results suggesting that the algorithm generates the set of well-distributed solutions approximating the Pareto front, even in the case for the problems with three objective functions.